\documentclass[twocolumn]{cifa}
\usepackage[latin1]{inputenc}
\usepackage[T1]{fontenc}
\usepackage[francais]{babel}
\usepackage{graphicx}
\usepackage{psfrag}

\usepackage{epsfig}
\usepackage[hang]{caption2}
\usepackage[normalsize,it,hang]{subfigure}
\usepackage{amsmath}
\usepackage{array}
\usepackage{theorem}
\usepackage{amsmath,amssymb}
\usepackage{amsfonts}
\usepackage{amssymb}
\usepackage{subfigure}
\usepackage{amssymb,amsbsy,amsmath,amsfonts,amssymb,amscd}

\newtheorem{remarque}{\it Remarque\/}

\newcommand{\Vect}[1]{ \mathbf{#1} }

\begin{document}

\title{Vers une commande multivariable \\ sans mod\`ele}
\author{\vskip 1em Michel \textsc{FLIESS}\textsuperscript{1}, C\'edric
\textsc{JOIN}\textsuperscript{2}, Mamadou
\textsc{MBOUP}\textsuperscript{3}, Hebertt \textsc{SIRA-RAM\'IREZ}\textsuperscript{4} \\
{\small \vskip1em \textsuperscript{1}Projet ALIEN, INRIA Futurs \&
\'Equipe MAX, LIX (UMR-CNRS
7161)\\
\'Ecole polytechnique, 91128 Palaiseau, France\\
\textsuperscript{2}Projet ALIEN, INRIA Futurs \& CRAN (UMR-CNRS
7039)
\\ Universit\'e Henri Poincar\'e (Nancy
I), BP 239, 54506 Vand\oe{}uvre-l\`es-Nancy, France\\
\textsuperscript{3}{Projet ALIEN, INRIA Futurs \& UFR de
Math\'ematiques et Informatique\\ Universit\'e Ren\'e-Descartes, 45
rue des Saints-P\`eres, 75270 Paris cedex 06, France} \\
\textsuperscript{4}{Departamento de Ingenier\'{\i}a El\'ectrica,
Secci\'on de Mecatr\'onica, Cinvestav-IPN \\ Aven. IPN No. 2508,
Colonia San Pedro Zacatenco, AP 14740, 07300 M\'exico, D.F.,
Mexique}}
{\small \vskip 1em { \texttt{Michel.Fliess@polytechnique.edu,
Cedric.Join@cran.uhp-nancy.fr \\
Mamadou.Mboup@math-info.univ-paris5.fr, hsira@cinvestav.mx}}}
} \maketitle

\begin{abstract}
On propose une commande de syst\`emes multivariables, de dimension
finie, lin\'eaires ou non, sans en conna\^{\i}tre le mod\`ele
math\'ematique. Nos deux outils essentiels sont l'alg\`{e}bre
diff\'{e}rentielle, et une estimation des d\'eriv\'ees de signaux
bruit\'es, r\'ecemment mise au point. Les simulations num\'eriques
de deux exemples, l'un lin\'{e}aire, l'autre non, valident notre
d\'emarche.
\end{abstract}

\begin{keywords}
Syst\`emes lin\'eaires multivariables, syst\`emes non lin\'eaires
multivariables, identification bo\^{\i}te noire, estimation, d\'{e}riv\'{e}es de
signaux bruit\'{e}s, alg\`ebre diff\'erentielle, calcul op\'{e}rationnel.
\end{keywords}

\section{Introduction}
Cette communication pr\'esente une g\'en\'eralisation multivariable
de \cite{kang}, o\`u \'etaient jet\'ees les bases d'une commande
sans mod\`ele pour syst\`emes monovariables de dimension finie,
lin\'eaires ou non. Est-il besoin de rappeler ({\it cf.}
\cite{ljung}) la difficult\'e redoutable d'obtenir un mod\`ele
math\'ematique fiable dans bien des situations quelque peu
complexes? Ainsi s'explique l'incroyable popularit\'e industrielle
des correcteurs PID, en d\'epit de r\'eglages souvent malais\'es
({\it cf.} \cite{pid,od}). Notre d\'emarche se distingue par sa
m\'ethodologie et sa \og philosophie \fg ~des identifications de
type \og bo\^{\i}te noire \fg, telles qu'on les trouve dans la
litt\'{e}rature (voir, par exemple, \cite{aut1} et \cite{aut2}). Elle
repose sur une nouvelle approche\footnote{Cette approche, n\'ee en
\cite{esaim} \`a propos d'identification param\'etrique lin\'eaire
en boucle ferm\'ee, a \'et\'e \'{e}tendue au signal \cite{mexi}. Elle
conduit \`{a} un changement radical de paradigme ({\it cf.}
\cite{bruit,jesa}). Voir, par exemple, \cite{b,linz} pour des
applications concr\`{e}tes.}, qui a permis l'estimation en temps r\'eel
des d\'eriv\'ees de signaux bruit\'es \cite{nolcos}\footnote{Voir
\cite{im,cras,gretsi,zeitz,nolcos,g,reger} pour les d\'{e}j\`a
nombreuses applications en automatique non lin\'eaire et traitement
du signal.}.

Nous substituons aux mod\`eles math\'ematiques d\'ecrivant les
machines dans une plage de fonctionnement aussi large que possible
des \'equations diff\'erentielles \og ph\'enom\'enologiques \fg,
valides sur un court laps de temps, et actualis\'ees pas \`a pas.
Nous \'ecrivons un syst\`eme multivariable, de dimension finie, \`a
$m$ entr\'ees ${\Vect{u}} = (u_1, \dots, u_m)$ et $p$ sorties
${\mathbf y} = (y_1, \dots, y_p)$, sous la forme {\small
\begin{equation} \label{mod} \begin{array}{l}y^{(n_1)}_1 = F_1 +
\alpha_{1,1} u_1 + \dots + \alpha_{1,m} u_m + \beta_1
\\ \ldots \\
y^{(n_p)}_p = F_p + \alpha_{p,1} u_1 + \dots + \alpha_{p,m} u_m +
\beta_p
\end{array}
\end{equation}}
o\`u
\begin{itemize}
\item $n_j \geq 1$, $j = 1, \dots, p$,  et, le plus souvent, $n_j = 1$, ou $2$;
\item $\alpha_{j,i}, \beta_{j} \in \mathbb{R}$, $i = 1, \dots, m$, $j = 1, \dots, p$, sont des param\`etres constants {\it non
physiques}, choisis par le praticien;
\item les $F_j$ sont d\'etermin\'es gr\^ace \`a la connaissance de $y^{(n_j)}_j$, $u_i$,
$\alpha_{j,i}$, $\beta_j$, $i = 1, \dots, m$, $j = 1, \dots, p$.
\item Si $p \gneqq m$, on garde seulement $m$ sorties de mani\`ere \`a
obtenir un syst\`eme carr\'e, inversible.
\end{itemize}
Le comportement d\'esir\'e s'obtient par correcteurs de type {\em
proportionnel int\'egral g\'en\'eralis\'e (GPI)} \cite{gpi} autour
d'une trajectoire de r\'ef\'erence. Dans les deux exemples
ci-dessous, on utilise, comme en \cite{kang}, un PID ou un PI.

Le {\S} \ref{io-eq} rappelle comment obtenir, gr\^ace \`a l'alg\`ebre
diff\'erentielle, les \'equations diff\'erentielles entr\'ee-sortie
d'un syst\`{e}me non lin\'{e}aire. Le {\S} \ref{estder} sur l'estimation des
d\'eriv\'ees d'un signal bruit\'e inclut certains des progr\`es
r\'ecents sur la mise en {\oe}uvre num\'erique. On \'evoque au {\S}
\ref{bbi} les principes essentiels de notre identification bo\^{\i}te
noire. Le {\S} \ref{exem} contient deux exemples, l'un lin\'eaire,
l'autre non, et leurs simulations num\'eriques. Une br\`eve
conclusion \'evoque quelques perspectives futures.

\section{Rappels sur les syst\`emes non lin\'eaires}\label{io-eq}
\subsection{Corps diff\'erentiels}
Un {\em corps diff\'erentiel}\footnote{Voir \cite{cl,kolchin} pour
plus de details et, en particulier, \cite{cl} pour des rappels sur
les corps usuels, c'est-\`a-dire non diff\'{e}rentiels. Tous les corps
consid\'er\'es ici sont de caract\'eristique nulle.} $\mathfrak{K}$
est un corps commutatif, muni d'une {\em d\'erivation}
$\frac{d}{dt}$, c'est-\`a-dire une application $\mathfrak{K}
\rightarrow \mathfrak{K}$ telle que, $\forall ~ a, b \in
\mathfrak{K}$,
{\small \begin{itemize}
\item $\frac{d}{dt}(a + b) = \dot{a} + \dot{b}$,
\item $\frac{d}{dt}(ab) = \dot{a}b + a\dot{b}$.
\end{itemize}}
Une {\em constante} $c \in \mathfrak{K}$ est un \'el\'ement tel que
$\dot{c} = 0$. L'ensemble des constantes est le {\em sous-corps des
constantes}.

Une {\em extension} de corps diff\'erentiels $\mathfrak{L} /
\mathfrak{K}$ consiste en la donn\'ee de deux corps diff\'erentiels
$\mathfrak{K}$, $\mathfrak{L}$, telles que:
\begin{itemize}
\item ${\mathfrak K} \subseteq {\mathfrak L}$,
\item la d\'erivation de $\mathfrak{K}$ est la restriction
\`a $\mathfrak{K}$ de celle de $\mathfrak{L}$.
\end{itemize}
Notons ${\mathfrak K}\langle S \rangle$, $S \subset {\mathfrak L}$,
le sous-corps diff\'erentiel de $\mathfrak{L}$ engendr\'e par
$\mathfrak{K}$ et $S$. Supposons $\mathfrak{L} / \mathfrak{K}$
finiment engendr\'e, c'est-\`a-dire ${\mathfrak L} = {\mathfrak
K}\langle S \rangle$, o\`u $S$ est fini. Un \'el\'ement $\xi \in
\mathfrak{L}$ est dit {\em diff\'erentiellement alg\'ebrique} par
rapport \`a $\mathfrak{K}$ si, et seulement si, il satisfait une
\'equation diff\'erentielle alg\'ebrique $P(\xi, \dots, \xi^{(n)}) =
0$, o\`u $P$ est un polyn\^ome sur $\mathfrak{K}$ en $n + 1$
ind\'etermin\'ees. L'extension $\mathfrak{L} / \mathfrak{K}$ est
dite {\em diff\'erentiellement alg\'ebrique} si, et seulement si,
tout \'el\'ement de $\mathfrak{L}$ de diff\'erentiellement
alg\'ebrique par rapport \`a $\mathfrak{K}$. Le r\'esultat suivant
est important: $\mathfrak{L} / \mathfrak{K}$ est
diff\'erentiellement alg\'ebrique si, et seulement si, son degr\'e
de transcendance est fini.

Un \'el\'ement de $\mathfrak{L}$ non diff\'erentiellement
alg\'ebrique par rapport \`a $\mathfrak{K}$ est dit {\em
diff\'erentiellement transcendant} par rapport \`a $\mathfrak{K}$.
Une extension $\mathfrak{L} / \mathfrak{K}$ non diff\'erentiellement
alg\'ebrique est dit {\em diff\'erentiellement transcendante}. Un
ensemble $\{ \xi_\iota \in \mathfrak{L} \mid \iota \in I \}$ est dit
{\em diff\'erentiellement alg\'ebriquement ind\'ependant} par
rapport \`a $\mathfrak{K}$ si, et seulement si, aucune relation
diff\'erentielle non triviale par rapport \`a $\mathfrak{K}$
n'existe: $Q( \xi_\iota^{(\nu_\iota)} ) = 0$, o\`u $Q$ est un
polyn\^ome sur $\mathfrak{K}$, implique $Q \equiv 0$. Deux ensembles
maximaux d'\'el\'ements diff\'erentiellement alg\'ebriquement
ind\'ependants ont m\^eme cardinalit\'e, c'est-\`a-dire m\^eme
nombre d'\'el\'ements: c'est le {\em degr\'e de transcendance
diff\'erentielle} de l'extension $\mathfrak{L} / \mathfrak{K}$. Un
tel ensemble est une {\em base de transcendance diff\'erentielle}.
Enfin, ${\mathfrak L} / {\mathfrak K}$ est diff\'erentiellement
alg\'ebrique si, et seulement si, son degr\'e de transcendance
diff\'erentielle est nulle.

\subsection{Syst\`emes non lin\'eaires}

Donnons-nous un corps diff\'erentiel de base $k$. Un {\em
syst\`eme}\footnote{Pour plus de d\'etails, voir
\cite{delaleau,flmr,rudolph,hsr}. Rappelons que l'on emploie aussi
l'alg\`{e}bre diff\'{e}rentielle dans diverses questions relatives \`{a}
l'identification et l'observabilit\'{e} (voir, par exemple,
\cite{ljung0}).} est une extension diff\'erentiellement
transcendante de $K/k$, finiment engendr\'ee. Soit $m$ son degr\'e
de transcendance diff\'erentielle. Un ensemble de {\em commandes
(independantes)} ${\Vect{u}} = (u_1, \dots, u_m)$ est une base de
transcendance diff\'erentielle de $K/k$. L'extension $K / k \langle
\Vect{u} \rangle$ est donc diff\'erentiellement alg\'ebrique. Un
ensemble de {\em sorties} ${\Vect{y}} = (y_1, \dots, y_p)$ est un
sous-ensemble de $K$.

Soit ${\Vect{x}} = (x_1, \dots, x_n)$ une base de transcendance de
$K / k \langle \Vect{u} \rangle$, de degr\'e de transcendance $n$.
Il en d\'ecoule la repr\'esentation d'\'etat g\'en\'eralis\'ee:
{\small $$
\begin{array}{l}
A_\iota(\dot{x}_\iota, {\Vect{x}}, {\Vect{u}}, \dots, {\Vect{u}}^{(\alpha)}) = 0 \\
B_\kappa(y_\kappa, {\Vect{x}}, {\Vect{u}}, \dots,
{\Vect{u}}^{(\beta)}) = 0
\end{array}
$$}

\noindent o\`{u} $A_\iota$, $\iota = 1, \dots, n$, $B_\kappa$, $\kappa =
1, \dots, p$, sont des polyn\^omes sur $k$.

La repr\'esentation entr\'ee-sortie suivante r\'esulte du fait que
$y_1, \dots, y_p$ sont diff\'erentiellement alg\'ebriques par
rapport \`a $k \langle {\mathbf u} \rangle$:
{\small
\begin{equation} \label{io} \Phi_j({\mathbf y}, \dots, {\mathbf
y}^{(\bar{N}_j)}, {\mathbf u}, \dots, {\mathbf u}^{(\bar{M}_j)} ) =
0
\end{equation}}
o\`u $\Phi_j$, $j = 1, \dots, p$, est un polyn\^ome sur $k$.

Venons-en \`a l'inversibilit\'e entr\'ee-sortie:
\begin{itemize}\item Le syst\`eme est dit {\em inversible \`a gauche}
si, et seulement si, l'extension $k \langle {\mathbf u}, {\mathbf y}
\rangle / k \langle {\mathbf y} \rangle $ est diff\'erentiellement
alg\'ebrique. C'est dire que l'on peut r\'ecup\'erer l'entr\'ee \`a
partir de la sortie gr\^ace \`a un syst\`eme diff\'erentiel. Alors,
$m \leq p$.
\item Il est dit {\em inversible \`a droite} si, et seulement si, le
degr\'e de transcendance diff\'erentielle de $k \langle {\mathbf y}
\rangle / k $ vaut $p$. C'est dire que les sorties sont
diff\'erentiellement alg\'ebriquement ind\'ependantes par rapport
\`a $k$. Alors, $p \leq m$.
\end{itemize}
Le syst\`eme est dit {\em carr\'e} si, et seulement si, $m = p$.
Alors, inversibilit\'es \`a gauche et \`a droite co\"{\i}ncident. On dit,
si ces propri\'et\'es sont v\'erifi\'ees, que le syst\`eme est {\em
inversible}.

{\small \begin{remarque} Supposons notre syst\`eme inversible \`a
gauche. Le comportement qualitatif de (\ref{io}), consid\'er\'e
comme syst\`eme d'\'equations diff\'erentielles en ${\mathbf u}$,
${\mathbf y}$ \'etant donn\'e, permet de d\'efinir le d\'ephasage
non minimal (voir aussi \cite{isidori2}).
\end{remarque}}

\section{Estimation des d\'eriv\'ees temporelles}\label{estder}
\`A la s\'erie convergente $x(t) = \sum_{n \geq 0} a_n \frac{t^n}{n
!}$, $a_n \in \mathbb{C}$, correspond la s\'erie op\'erationnelle
convergente\footnote{Voir \cite{miku1,miku2}. Renvoyons \`{a}
\cite{cras,gretsi} pour plus de d\'{e}tails.} $x = \sum_{n \geq 0}
\frac{a_n}{s^{n + 1}}$. Avec le d\'eveloppement de Taylor tronqu\'e
$x_N(t) \sum_{n=0}^N a_n \frac{t^n}{n !}$, l'\'equation
diff\'erentielle $\frac{d^{N+1}}{dt^{N+1}} = x_N = 0$ correspond
dans le domaine op\'erationnel \`a  {\small $$s^{N+1} x_N - s^{N}
x_N(0) - s^{N-1} \dot{x}_N(0)\ldots - x_N^{(N)}(0) = 0$$} Les
d\'eriv\'ees \`a  l'origine  $x_N^{(i)}(0)$ sont ainsi obtenues \`a
partir du syst\`eme d'\'equations lin\'eaires
{\small
\begin{multline}\label{eqdiff}
  s^{-\nu} \frac{d^m}{ds^m}\left\{x_N^{(N)}(0) + x_N^{(N-1)}(0) s
+ \ldots  + x_N(0) s^N\right\} =\\ s^{-\nu}
\frac{d^m}{ds^m}\left\{s^{N+1} x_N\right\}
\end{multline}}

\noindent $m = 0, \ldots, N$, $\nu \geqslant N+1$. Ce syst\`eme
\'etant triangulaire avec des \'el\'ements diagonaux non nuls, les
param\`etres $x_N^{(i)}(0)$, et, par cons\'equent, les coefficients
$a_0, \dots, a_N$ sont {\em lin\'eairement identifiables}
\cite{cras,esaim}. Rempla\c{c}ons $x_N$ par $x$ dans \eqref{eqdiff} : on
obtient ainsi l'{\em estim\'ee op\'erationnelle}
$[x^{(i)}(0)]_{e_N}$ de $x^{(i)}(0)$.

Pour le passage au num\'erique, il suffit, selon les r\`egles
usuelles du calcul op\'erationnel ({\it cf.}
\cite{miku1,miku2,pol}), de remplacer en (\ref{eqdiff})
\begin{itemize}
\item $\frac{c}{s^\alpha}$, $\alpha \geq 1$, $c \in \mathbb{C}$, par
$c \frac{t^{\alpha - 1}}{(\alpha - 1) !}$, $t \geq 0$;
\item $\frac{1}{s^\alpha}\frac{d^n x}{ds^n}$ par l'int\'egrale it\'er\'ee
d'ordre $\alpha$ {\small \begin{multline}
\int_0^{t}\!\!\int_0^{t_{\alpha-1}}\cdots\int_0^{t_1}
(-1)^n \tau^n x(\tau) d t_{\alpha-1}\cdots d t_1 d\tau =\\
\frac{(-1)^n}{(\alpha - 1) !}\int_0^t (t-\tau)^{\alpha - 1} \tau^n
x(\tau) d\tau
\end{multline}}
\end{itemize}
Notons, $[x^{(i)}(0)]_{e_N}(t)$ l'{\em estim\'ee num\'erique}  ainsi
obtenue de $x^{(i)}(0)$, pour un temps d'estimation $t$. La mise en
{\oe}uvre repose sur le r\'esultat suivant: {\small $$ \lim_{t
\downarrow 0} [x^{(i)}(0)]_{e_N}(t) = \lim_{N \to +\infty}
[x^{(i)}(0)]_{e_N}(t) = x^{(i)}(0)
$$}
{\small \begin{remarque} Les it\'erations des int\'egrales
produisent une moyennisation, donc un filtrage passe-bas, qui permet
d'att\'enuer les bruits (voir \cite{bruit}).
\end{remarque}
\begin{remarque}
La fen\^etre temporelle d'estimation peut \^etre choisie tr\`es
petite, ce qui permet une impl\'ementation en temps r\'eel.
\end{remarque}}

\section{Proc\'edures d'identification bo\^{\i}te noire}\label{bbi}
On requiert les propri\'et\'es suivantes pour \'eviter, notamment,
toute boucle alg\'ebrique:
\begin{enumerate}
\item On suppose le syst\`eme inversible \`a gauche. Si le nombre de
sorties est strictement sup\'erieur \`a celui des entr\'ees,
c'est-\`a-dire $p \gneqq m$, on choisit $m$ sorties pour obtenir un
syst\`eme carr\'e inversible, et construire (\ref{mod}).

\item  L'ordre de d\'erivation $n_j$ en (\ref{mod}) est reli\'e \`a
(\ref{io}) par $n_j \leq \bar{N}_j$ et, plus pr\'ecis\'ement, par
$\frac{\partial \Phi_j}{\partial y^{(n_j)}} \not\equiv 0$.

\item La valeur num\'erique de $F_j$, \'egale \`a $y^{(n_j)}_j - \alpha_{j,1} u_1
- \dots - \alpha_{j,m} u_m - \beta_j$ est obtenue gr\^ace \`a la
discr\'etisation {\small \begin{equation*}\label{Fk} F_{j}(\kappa) =
[y^{(n_j)}_j(\kappa)]_e - \sum_{i = 1}^m \alpha_{j,i} u_i(\kappa -
1) - \beta_j \end{equation*}} o\`u $[ \bullet (\kappa) ]_e$
d\'esigne l'estim\'ee au temps $\kappa$.
\end{enumerate}
Le praticien suit les \'etapes suivantes:
\begin{enumerate}
\item choix de $m$ sorties si $p > m$;
\item choix des $n_j$ en (\ref{mod});
\item choix de trajectoires de r\'ef\'erence pour les $y_j$, ainsi qu'il est
usuel en commande par platitude (voir, par exemple,
\cite{flmr,rudolph,hsr})\footnote{C'est donc une commande pr\'{e}dictive
sans mod\`{e}le. Renvoyons \`a \cite{mfrm} et \cite{delhag} pour les
avantages de la platitude en pr\'{e}dictif avec mod\`{e}le.};

\item choix des param\`etres $\alpha_{j,i}$ en (\ref{mod}) afin que que la grandeur
des commandes $\Vect{u}$ soient convenables;

\item choix des param\`etres $\beta_j$ en (\ref{mod}) diff\'erent de $0$ si $u$
n'appara\^{\i}t pas lin\'eairement en (\ref{io}), c'est-\`a-dire {\small
$\frac{\partial \Phi_j}{\partial u_j} ({\mathbf u} = 0) \equiv 0$}.
\end{enumerate}

{\small \begin{remarque} Avec des syst\`emes \`a d\'ephasage non
minimal, nos proc\'edures peuvent conduire \`a des valeurs
divergentes des $u_j$ pour $t$ grand, et, donc, \`a des valeurs
num\'eriquement inadmissibles des $F_j$.
\end{remarque}}


\section{Deux exemples}\label{exem}
On utilise les mod\`eles math\'ematiques ci-dessous pour les besoins
\'evidents des simulations num\'eriques.

\subsection{Syst\`eme lin\'eaire}\label{ex1}
Soit le syst\`eme lin\'eaire \`a deux entr\'ees et deux sorties,
avec p\^oles instables et large spectre:
{\small $$\begin{array}{ll}
y_1=&\frac{s^3}{(s+0.01)(s+0.1)(s-1)s}u_1\\
y_2=&\frac{s+1}{(s+0.003)(s-0.03)(s+0.3)(s+3)}u_1\\&+\frac{s^2}{(s+0.004)(s+0.04)(s-0.4)(s+4)}u_2
\end{array}$$}

\begin{figure*}[!ht]
\centering
%
{\subfigure[\footnotesize Bruits de sortie]{\rotatebox{-90}{\resizebox{!}{8.5cm}{\includegraphics{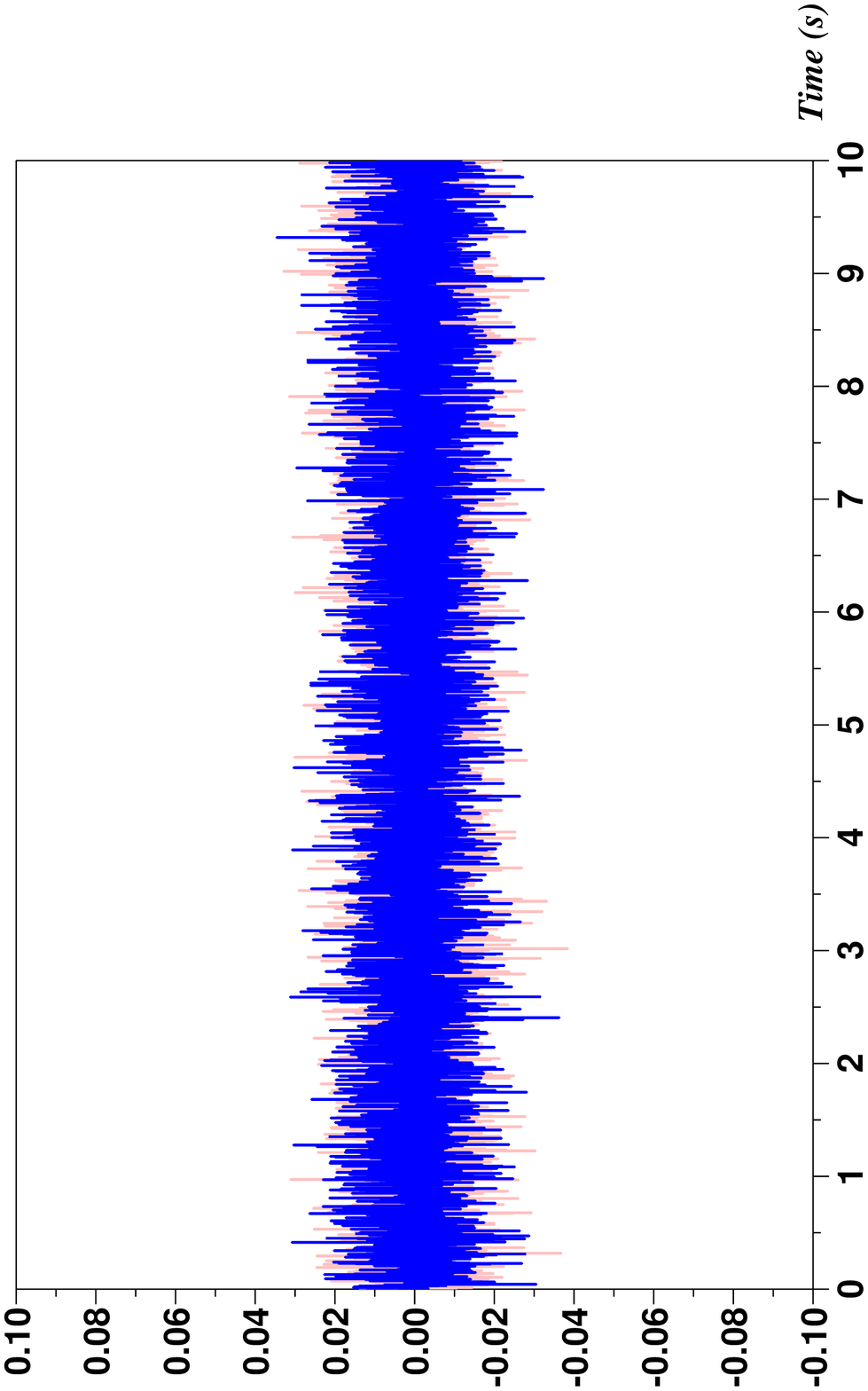}}}}}
{\subfigure[\footnotesize Estimation de $\dot y_1$]{
\rotatebox{-90}{\resizebox{!}{8.5cm}{\includegraphics{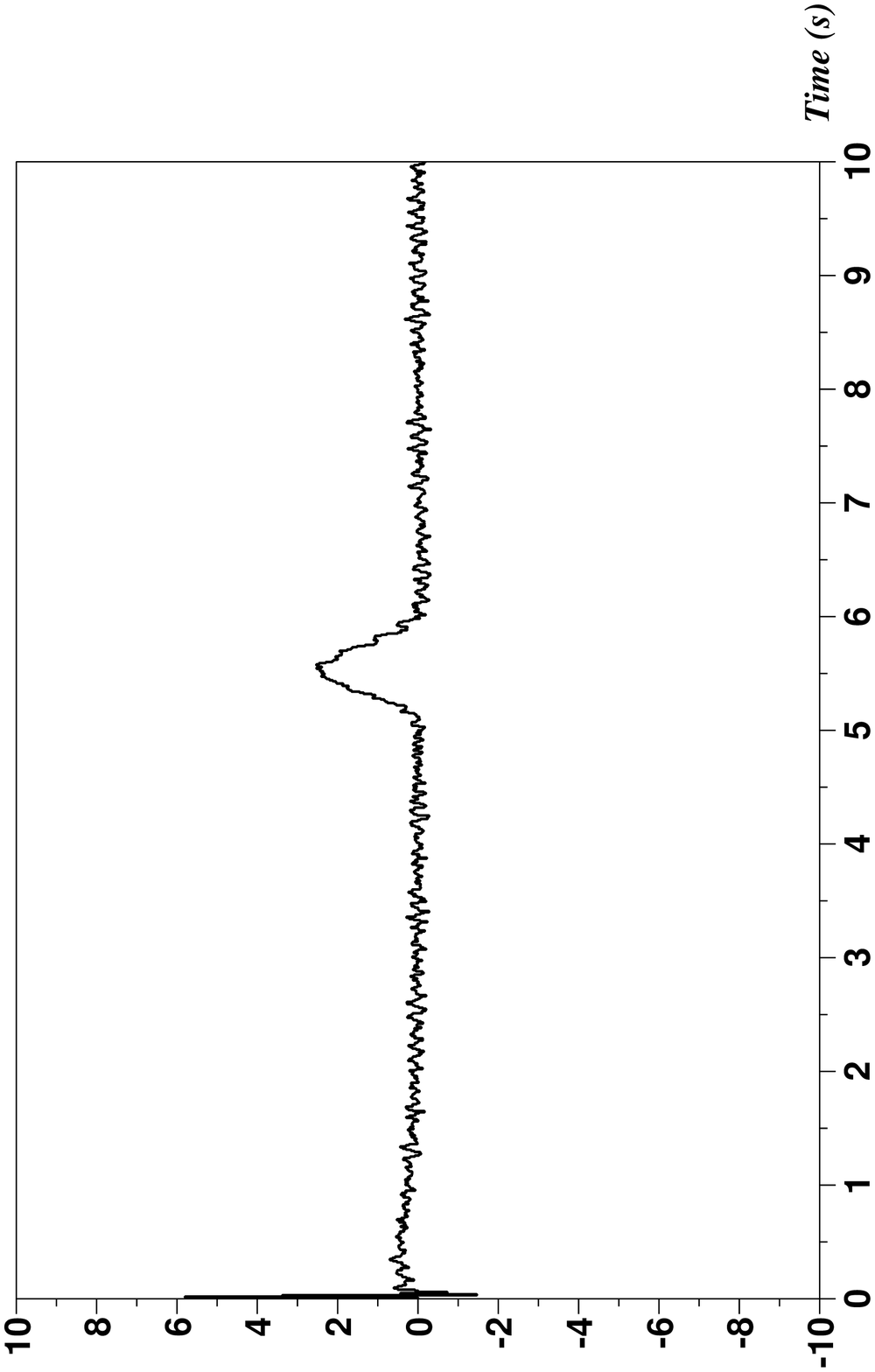}}}}}
{\subfigure[\footnotesize Estimation de $\dot y_2$]{
\rotatebox{-90}{\resizebox{!}{8.5cm}{\includegraphics{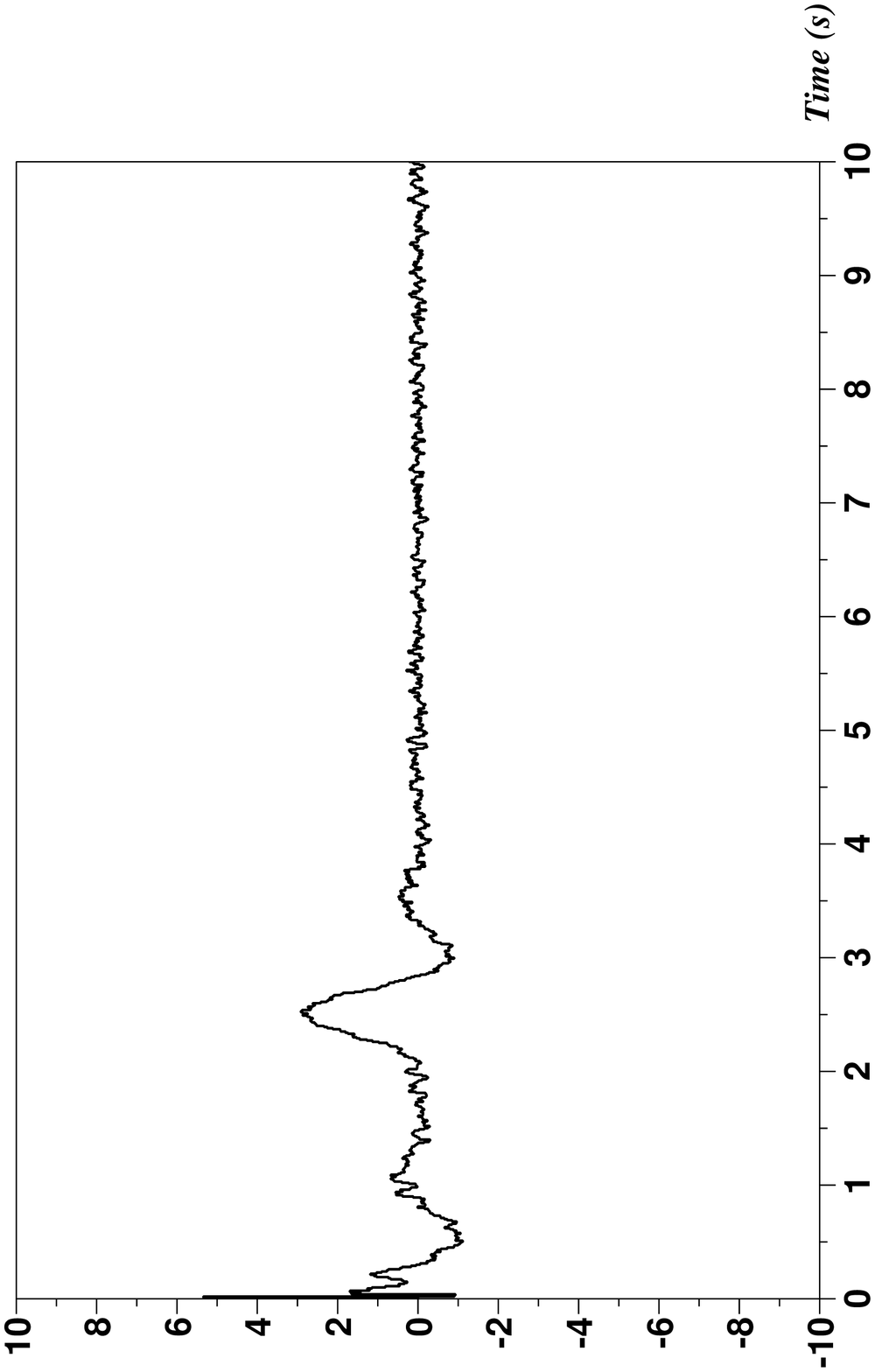}}}}}
%
{\subfigure[\footnotesize Estimation de $\ddot y_2$]{
\rotatebox{-90}{\resizebox{!}{8.5cm}{\includegraphics{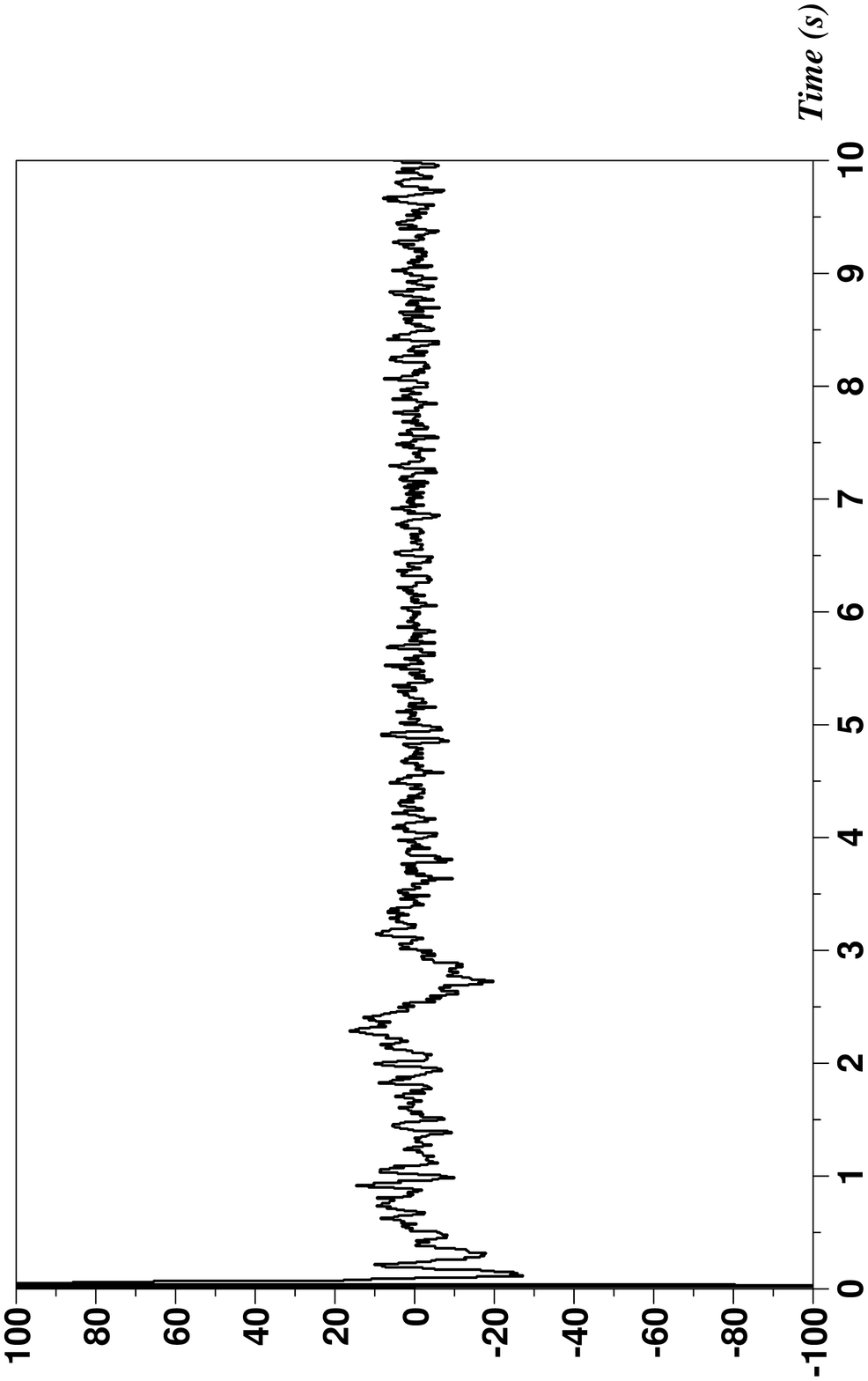}}}}}
{\subfigure[\footnotesize Estimation de $F_1$]{
\rotatebox{-90}{\resizebox{!}{8.5cm}{\includegraphics{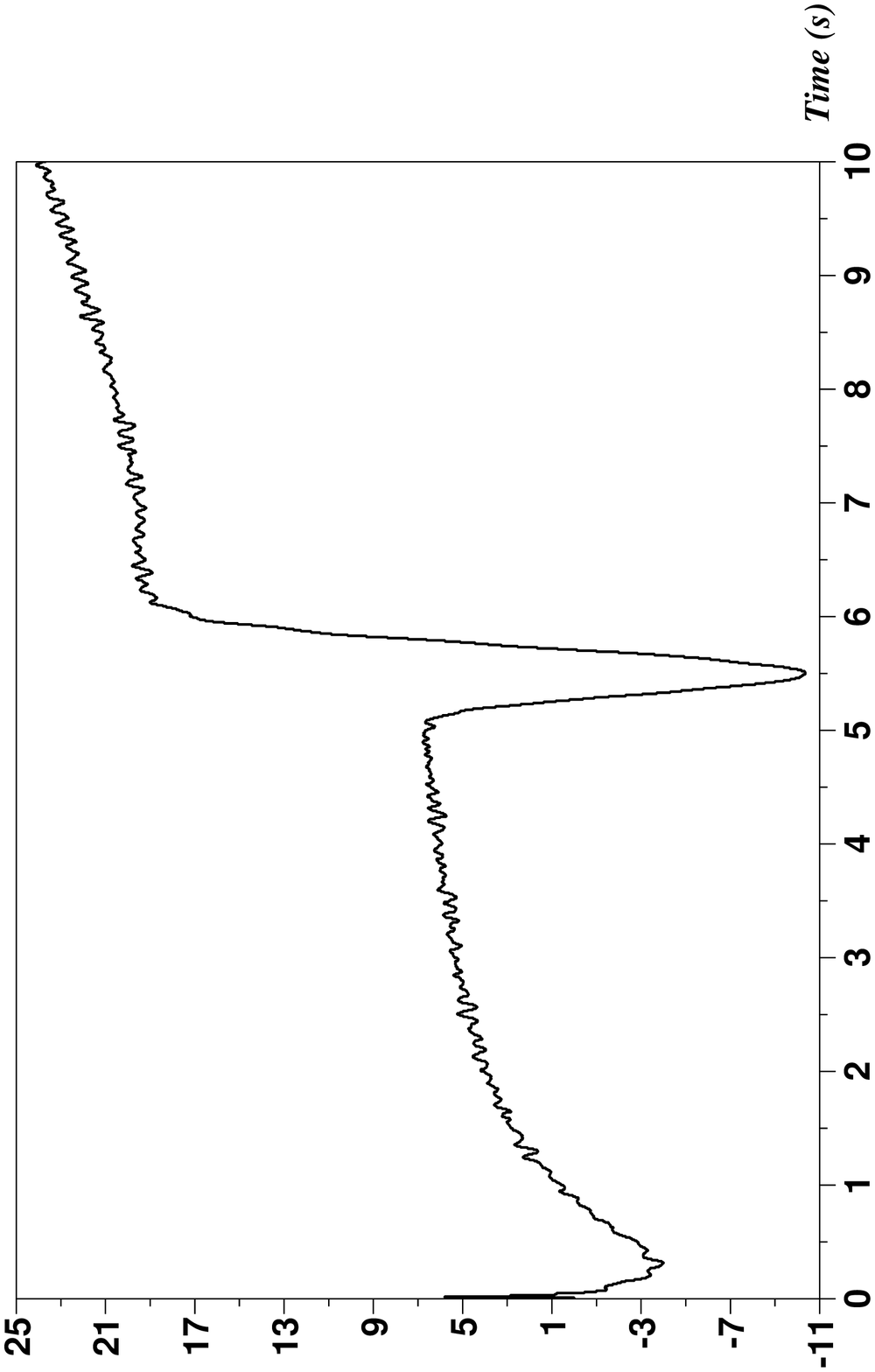}}}}}
{\subfigure[\footnotesize Estimation de $F_2$]{
\rotatebox{-90}{\resizebox{!}{8.5cm}{\includegraphics{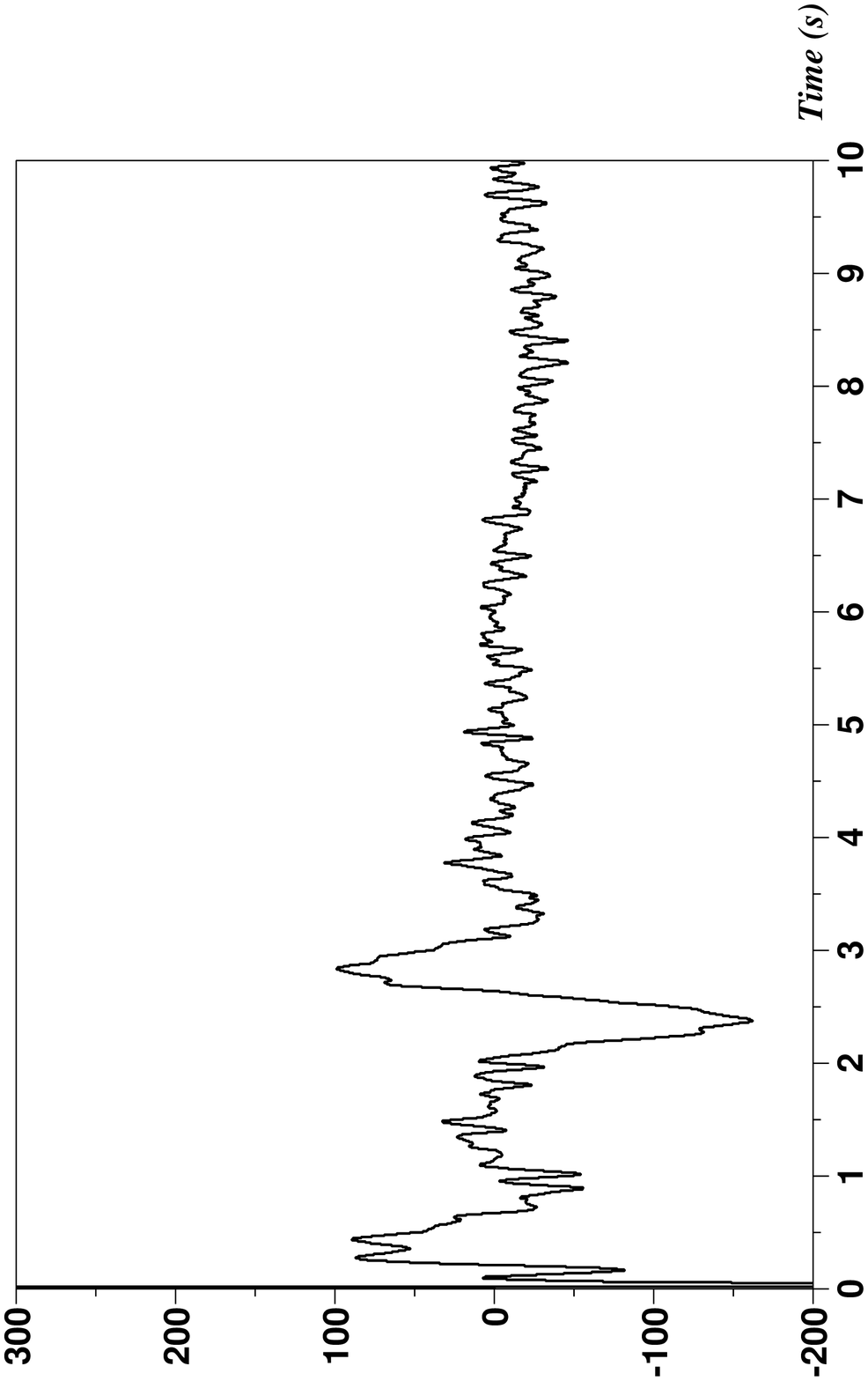}}}}}
{\subfigure[\footnotesize Estimation de $u_1$]{
\rotatebox{-90}{\resizebox{!}{8.5cm}{\includegraphics{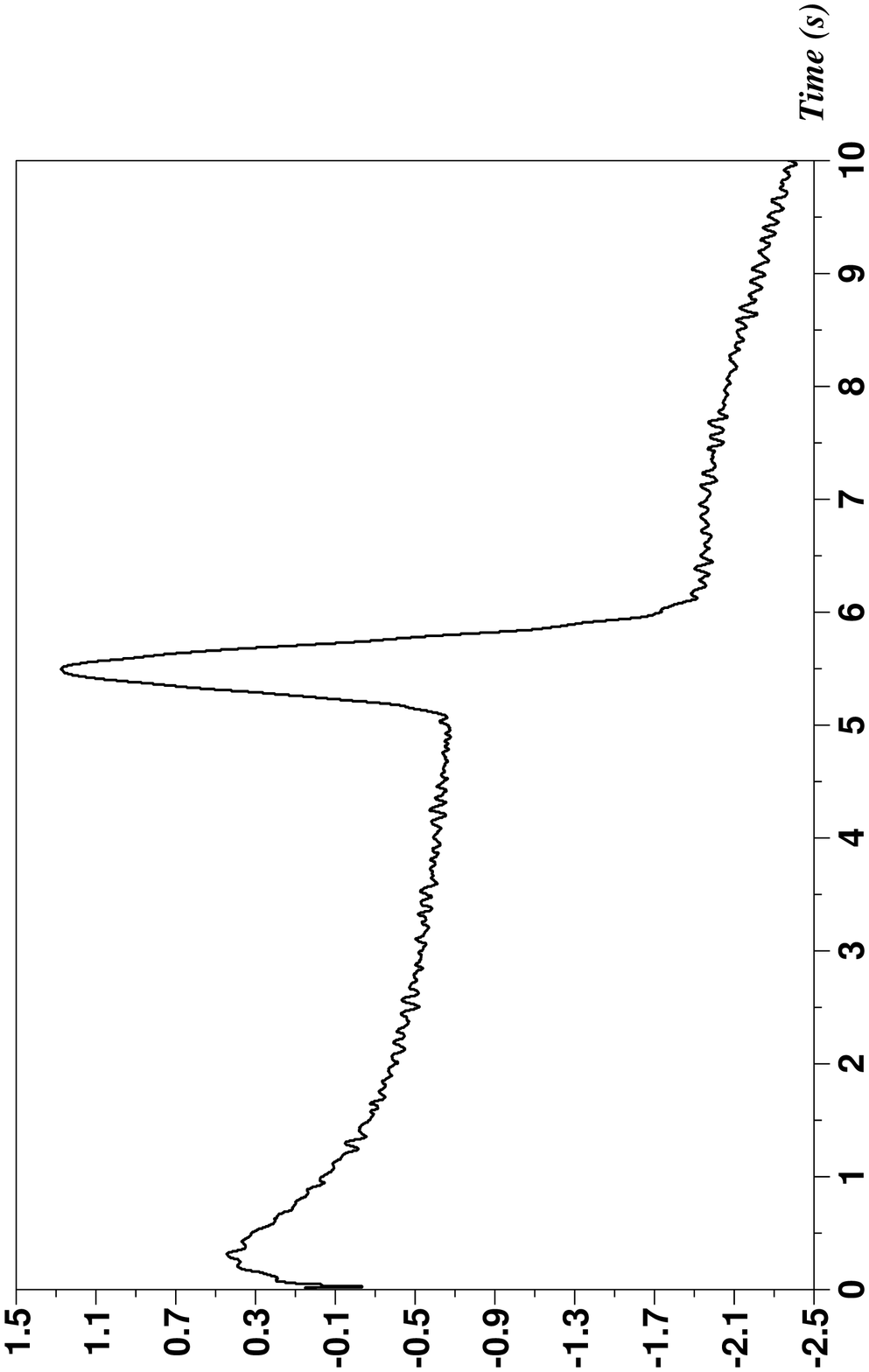}}}}}
{\subfigure[\footnotesize Estimation de $u_2$]{
\rotatebox{-90}{\resizebox{!}{8.5cm}{\includegraphics{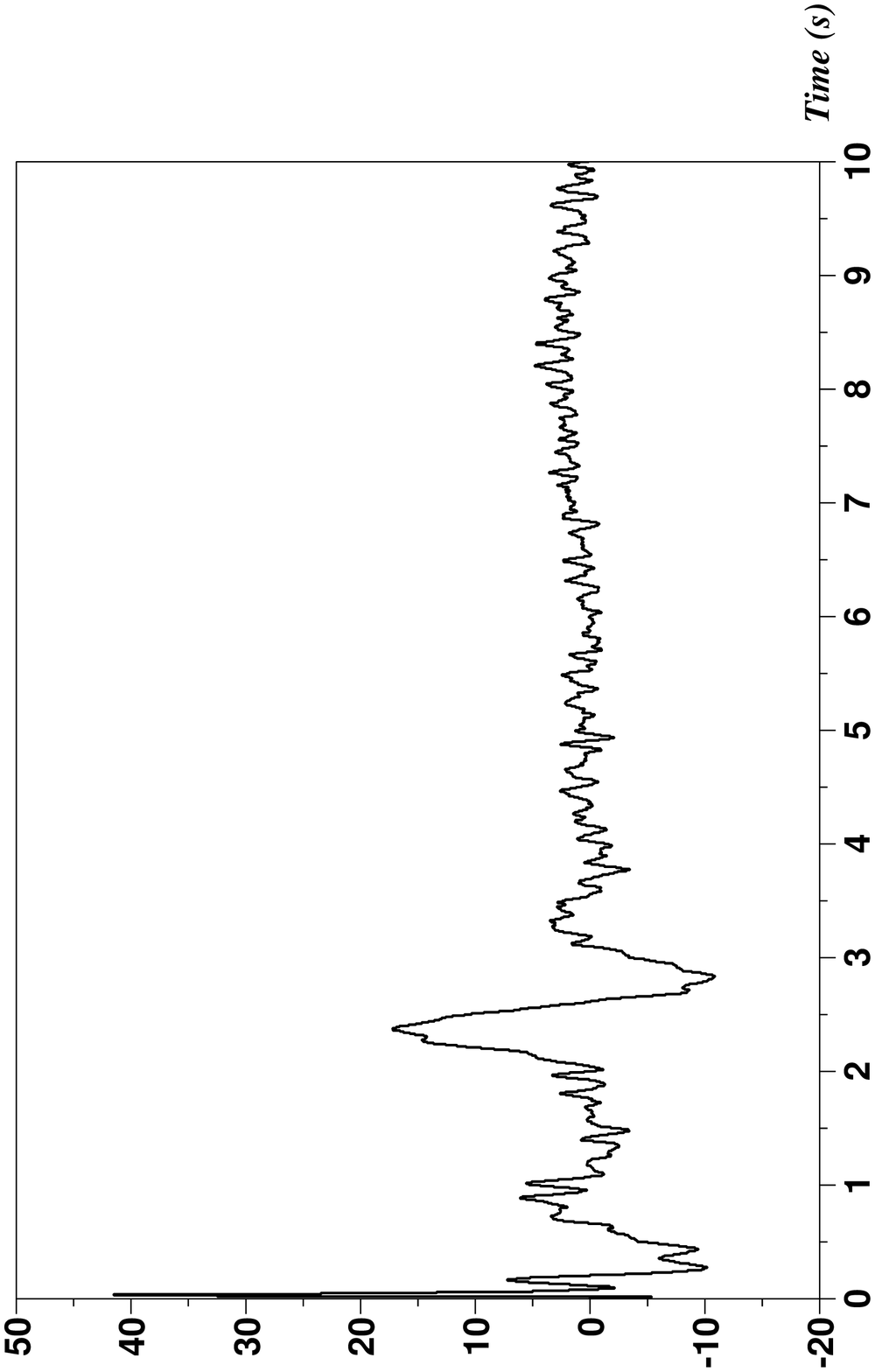}}}}}
 %
%
 \caption{Syst\`eme lin\'eaire\label{fig_Ml1}}
\end{figure*}

\noindent Apr\`es quelques essais, nous choisissons (\ref{mod}) sous
la forme d\'ecoupl\'ee: {\small \begin{equation*}\label{li}
\dot{y}_1=F_1 + 10u_1 \quad ~ \ddot{y}_2=F_2 + 10u_2
\end{equation*}}
La stabilisation autour d'une trajectoire de r\'ef\'erence est
assur\'ee par un r\'egulateur PID (voir \cite{kang}):
{\small
\begin{equation}\label{cor}
\begin{array}{l}
u_1 =\frac{1}{10}\left( \dot{y}^{\ast}_1 - F_1 + K_{P1} e_1 + K_{I1} \int e_1 +K_{D1} \dot{e}_1 \right)\\
u_2=\frac{1}{10}\left(\ddot{y}^{\ast}_2 -F_2 + K_{P2} e_2 +
K_{I2} \int e_2 +K_{D2} \dot{e}_2 \right)
\end{array}
\end{equation}}
o\`u \begin{itemize} \item {\small $K_{P1}=1$, $K_{I1}=K_{D1}=0$,
$K_{P2}=K_{I2}=50$, $K_{D2}=10$}; \item $y^{\ast}_1$, $y^{\ast}_2$
sont les trajectoires de r\'ef\'erence;
\item $e_1 = y^{\ast}_1 - y_1$, $e_2 = y^{\ast}_2 - y_2$.
\end{itemize} Le comportement en suivi de trajectoires, avec bruit de sortie
additif (loi normale $N(0,0.01)$), est bon\footnote{
C'est pourquoi nous esp\'erons, comme d\'ej\`a dit en \cite{kang},
que nos m\'ethodes pourraient fournir une alternative efficace \`a
celles sur la r\'eduction de mod\`eles (voir, par exemple,
\cite{antoulas}).}. Les estimations des signaux n\'ecessaires \`a la
synth\`ese de la commande sont pr\'esent\'ees dans les figures
\ref{fig_Ml1}-(b) \`a \ref{fig_Ml1}-(h). Notons la diff\'erence
d'\'echelle entre $u_1$ (figure \ref{fig_Ml1}-(g)) et $u_2$ (figure
\ref{fig_Ml1}-(h)), d'o\`u
l'in\'egalit\'e 
des param\`etres des correcteurs (\ref{cor}).

\begin{figure*}[!ht]
\centering {\subfigure[\footnotesize R\'ef\'erences (- -) et
sorties]{
\rotatebox{-90}{\resizebox{!}{8.5cm}{%
   \includegraphics{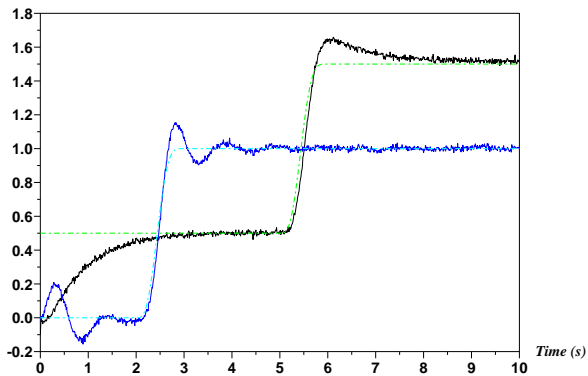}}}}}
{\subfigure[\footnotesize Sorties sans utiliser $F_1$ et $F_2$]{
\rotatebox{-90}{\resizebox{!}{8.5cm}{%
   \includegraphics{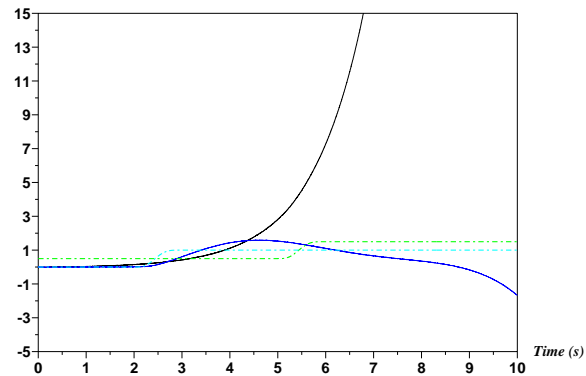}}}}}
 \caption{Syst\`eme lin\'eaire (suite)\label{fig_Ml2}}
\end{figure*}

La figure \ref{fig_Ml2}-(b) montre le comportement du syst\`eme en
appliquant une commande PID plus \og traditionnelle \fg,
c'est-\`a-dire en posant $F_1 \equiv F_2 \equiv 0$. La comparaison
des figures \ref{fig_Ml2}-(a) et \ref{fig_Ml2}-(b) est \'eloquente.

\begin{figure*}[!ht]
\centering

        \psfrag{p1}[cc][t]{\footnotesize pompe $1$}
        \psfrag{p2}[cc][t]{\footnotesize pompe $2$}
        \psfrag{q1}[l][cc]{\footnotesize$u_1$}
        \psfrag{q2}[r][cc]{\footnotesize$u_2$}
        \psfrag{c1}[r][cc]{\footnotesize cuve $1$}
        \psfrag{c2}[r][cc]{\footnotesize cuve $2$}
        \psfrag{c3}[r][cc]{\footnotesize cuve $3$}
        \psfrag{s}[l][cc]{$\tiny S$}
%
        \psfrag{q13}[cc][cc]{\tiny}
        \psfrag{q32}[cc][cc]{\tiny}
        \psfrag{q20}[cc][cc]{\tiny}
        \psfrag{m1}[l][cc]{\footnotesize$S_p,\mu_1$}
        \psfrag{m2}[l][cc]{\footnotesize$S_p,\mu_2$}
        \psfrag{m3}[l][cc]{\footnotesize$S_p,\mu_3$}
%
%
        \psfrag{l1}[r][cc]{\footnotesize$x_1$}
        \psfrag{l2}[r][cc]{\footnotesize$x_2$}
        \psfrag{l3}[r][cc]{\footnotesize$x_3$}
%
{\includegraphics[scale=0.28]{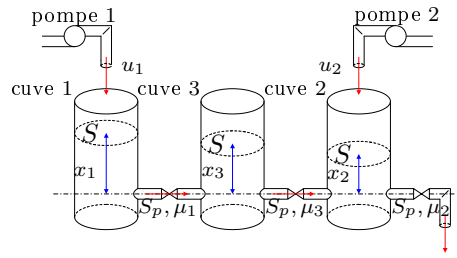}} \caption{Syst\`eme non
lin\'{e}aire des trois cuves \label{schem}}
\end{figure*}

\subsection{Syst\`eme non lin\'eaire des trois cuves}
Le syst\`eme des trois cuves de la figure \ref{schem}, tr\`es
populaire dans la communaut\'e du diagnostic\footnote{Renvoyons \`a
\cite{zeitz} pour plus de d\'etails et les r\'ef\'erences
bibliographiques. Que l'on nous permette de rappeler que
\cite{zeitz} propose, sans doute pour la premi\`ere fois, le
diagnostic, la commande et la reconfiguration d'un syst\`eme non
lin\'eaire \`a param\`etres incertains.}, v\'erifie les \'equations:
{\small $$
    \begin{cases}\begin{array}{ll}
    \dot{x}_1 =&-C_1\text{sign}(x_1 - x_3)\sqrt{|x_1 - x_3|}+u_1 /S\\
    \dot{x}_2 =&C_3\text{sign}(x_3 - x_2)\sqrt{|x_3 - x_2|}\\&-C_2\text{sign}(x_2)\sqrt{|x_2|}
    +u_2 /S\\
    \dot{x}_3 =&C_1\text{sign}(x_1 - x_3)\sqrt{|x_1 - x_3|}\\&-C_3\text{sign}(x_3 - x_2)\sqrt{|x_3 - x_2|}\\
    y_{1} =&x_1\\
    y_{2} =&x_2\\
    y_{3} =&x_3\\
    \end{array}\end{cases}
$$}
o\`{u} {\small $$
\begin{array}{l}
C_n = (1/S).\mu_n.S_p\sqrt{2g}, n = 1, 2, 3~;\\
S=0.0154~m~\text{(section des cuves)}~;\\
Sp=5.10^{-5}~m~\text{(section des tuyaux inter-cuves)}~;\\
g=9.81~m.s^{-2}~\text{(acc\'el\'eration de la pesanteur)}~;\\
\mu_1=\mu_3=0.5\text{,} ~ \mu_2=0.675 ~\text{(coefficients de
viscosit\'e)}.
\end{array}
$$}

\begin{figure*}[!ht]
\centering {\subfigure[\footnotesize R\'ef\'erences (- -) et
sorties]{
\rotatebox{-90}{\resizebox{!}{8.5cm}{%
   \includegraphics{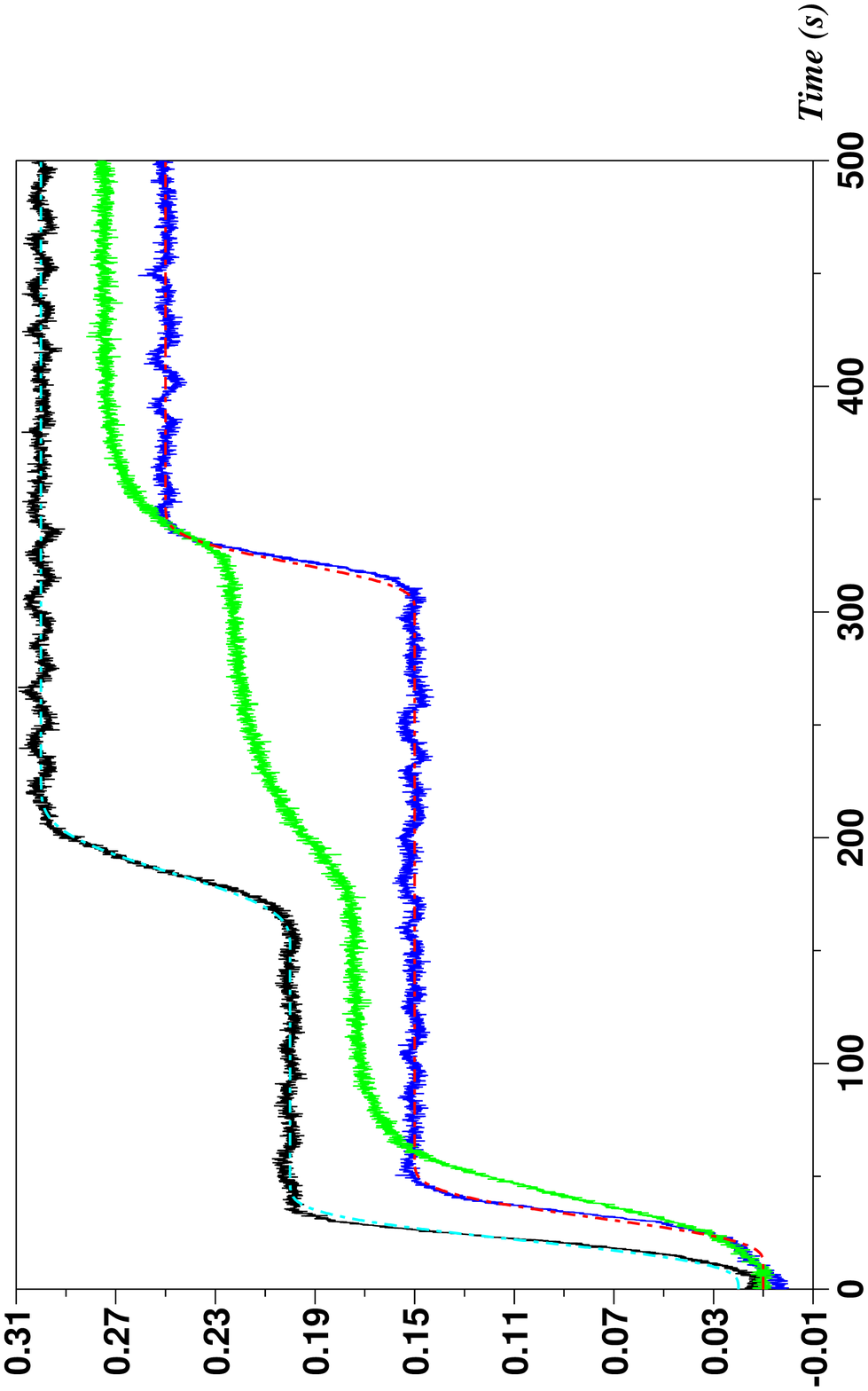}}}}}
{\subfigure[\footnotesize Estimation des d\'eriv\'ees $\dot y_1$ (-)
et $\dot y_2$ (- -), d\'{e}cal\'{e} de $0.02$]{
\rotatebox{-90}{\resizebox{!}{8.5cm}{\includegraphics{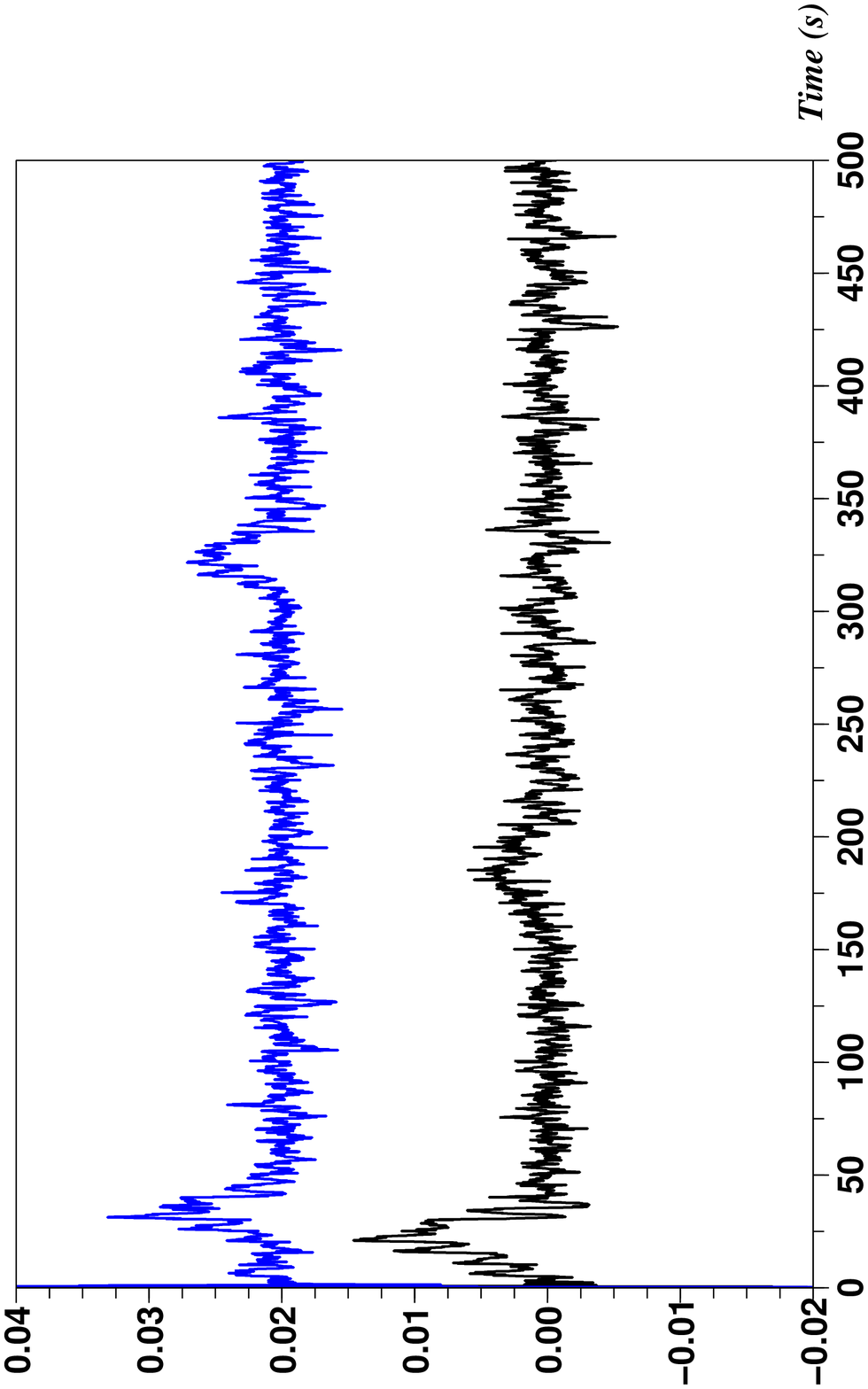}}}}}
{\subfigure[\footnotesize Commandes $u_1$ (-) et $u_2$ (- -)]{
\rotatebox{-90}{\resizebox{!}{8.5cm}{%
   \includegraphics{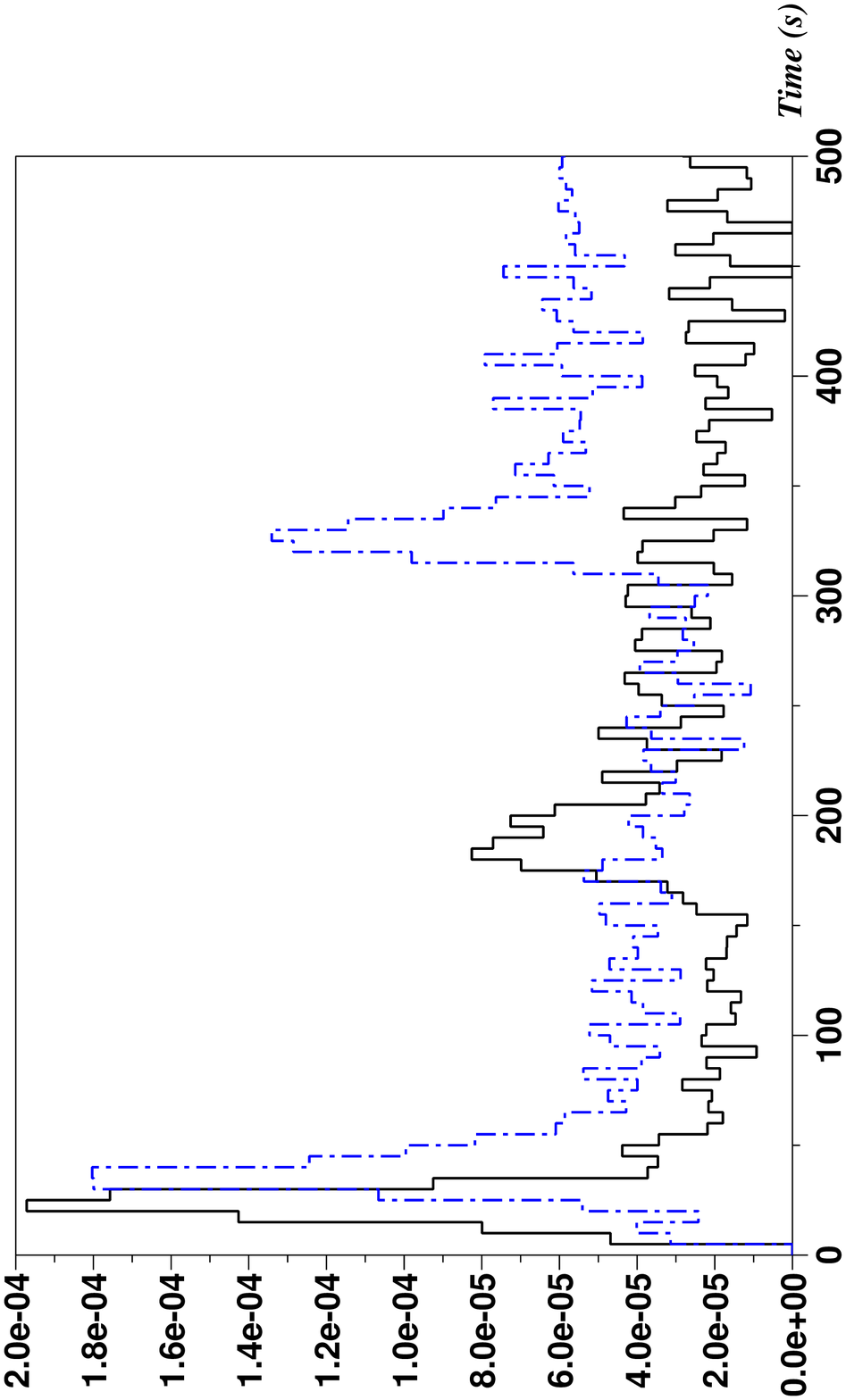}}}}}
{\subfigure[\footnotesize Gros plan sur le d\'ebruitage]{
\rotatebox{-90}{\resizebox{!}{8.5cm}{%
   \includegraphics{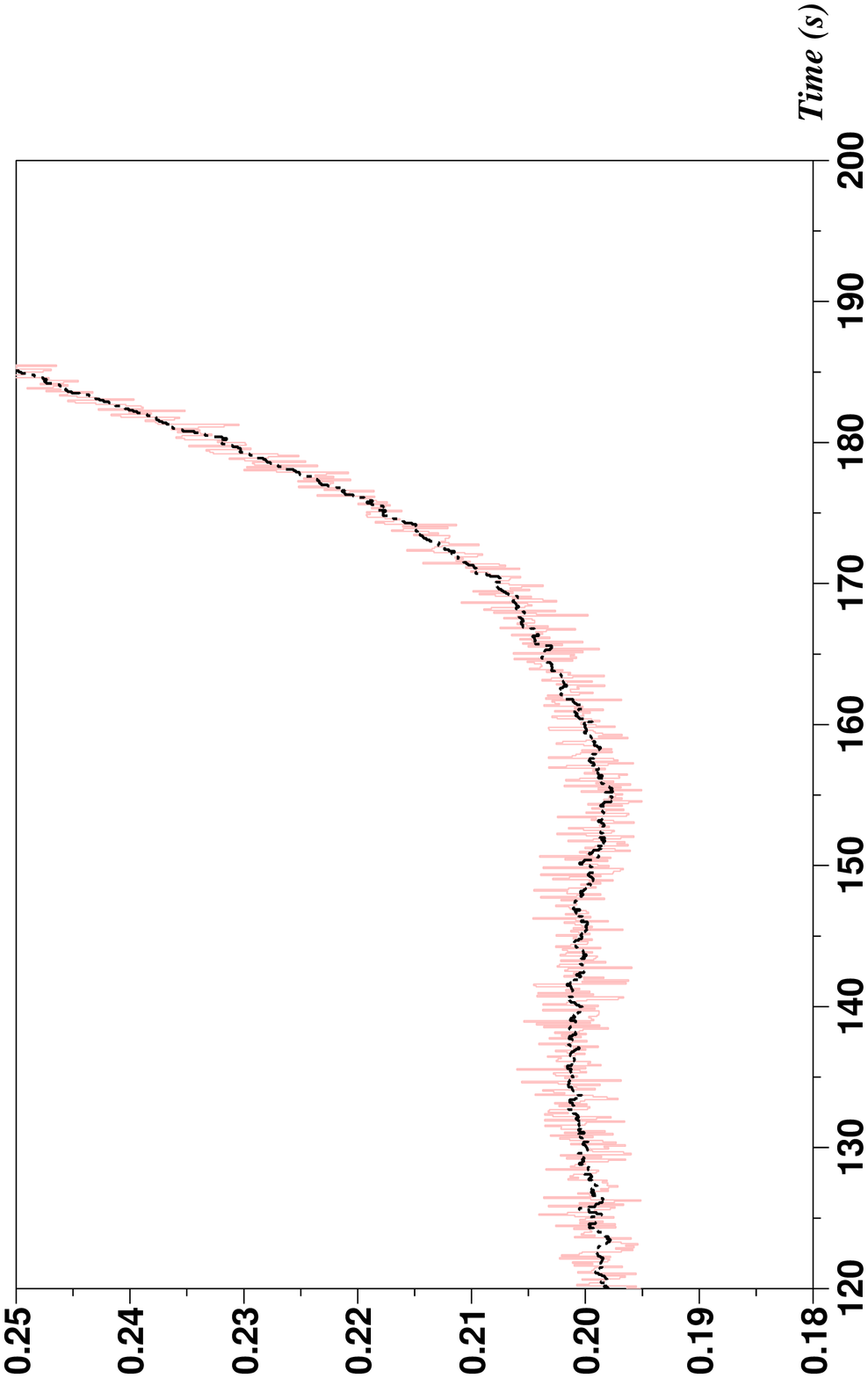}}}}}
 \caption{Syst\`eme non lin\'eaire des trois cuves\label{fig_3cuves}}
\end{figure*}

\noindent Selon les recommandations du {\S} \ref{bbi}, on construit
(\ref{mod}), d\'ecoupl\'e comme au {\S} \ref{ex1}: $\dot y_i=F_i+ 200
u_i$,  $i = 1, 2$. La figure \ref{fig_3cuves}-(a) fournit le suivi
de trajectoires. L'estimation des d\'eriv\'ees (figure
\ref{fig_3cuves}-(b)) poss\`ede un comportement remarquable en
d\'epit du bruit additif de mesure, de m\^eme caract\'eristique
qu'au {\S} \ref{ex1}. Les commandes nominales (figure
\ref{fig_3cuves}-(c)) sont assez proches de celles que nous aurions
calcul\'ees en utilisant la platitude (voir \cite{zeitz}). Elles
sont compl\'et\'ees par des correcteurs PI

{\small
\begin{equation*}\label{pid} u_i =
\frac{1}{200}\left(\dot{y}_{i}^\ast-F_i + 10 e_i + 2.10^{-2} \int
e_i \right) \quad i = 1, 2 \end{equation*}}

\noindent o\`u $y_{i}^\ast$ est la trajectoire de r\'ef\'erence,
$e_i = y_{i}^\ast - y_i$. Pour \'evaluer $e_i$ nous utilisons $y_i$
d\'ebruit\'e (voir figure \ref{fig_3cuves}-(d)) selon les techniques
du {\S} \ref{estder} (voir aussi \cite{cras,gretsi}).

\section{Conclusion}
Les r\`egles ardues d'identification du {\S} \ref{bbi} seront
pr\'ecis\'ees dans le futur\footnote{Elles ne peuvent \^etre toutes
d\'eduites de consid\'erations purement math\'ematiques.
L'exp\'erience pratique y joue bien entendu, comme pour les PID, un
r\^ole consid\'erable.}. On exposera bient\^ot des r\'{e}sultats
encourageants sur le d\'{e}phasage non minimal ainsi que ceux sur l'{\em
\'{e}galisation aveugle} ({\it cf.} \cite{eg}), qui est, en un certain
sens, le pendant en signal de la commande sans mod\`{e}le.

Une math\'{e}matisation \'{e}labor\'{e}e, comme celle pr\'{e}sent\'{e}e ici, afin
d'abandonner une mod\'{e}lisation aussi \og globale \fg ~que possible
dans une discipline empirique, comme l'automatique, semble
nouvelle\footnote{Voir, cependant, l'utilisation des dynamiques
lentes-rapides pour simplifier la mise en \'{e}quation ({\it cf.}
\cite{bo}).}. Il conviendrait d'en explorer les implications
\'{e}pist\'{e}mologiques, ne serait-ce que pour les liens entre physique,
complexit\'{e}, commande et r\'{e}solution temporelle ({\it cf.}
\cite{not}).


\end{document}